\documentclass[11pt,a4paper]{article}
\usepackage[utf8]{inputenc}
\usepackage[T1]{fontenc}
\usepackage[english]{babel}
\usepackage{geometry}
\usepackage{amsmath,amssymb,amsthm,mathtools}

\usepackage{booktabs}
\usepackage{cite}
\usepackage[colorlinks,citecolor=red,urlcolor=blue,linkcolor=blue]{hyperref}
\usepackage[capitalize,noabbrev]{cleveref}
\usepackage{doi}

\usepackage{verbatim}
\usepackage{comment}

\numberwithin{equation}{section}
\theoremstyle{plain}
\newtheorem{theorem}{Theorem}[section]
\newtheorem{proposition}[theorem]{Proposition}
\newtheorem{lemma}[theorem]{Lemma}
\newtheorem{corollary}[theorem]{Corollary}
\theoremstyle{remark}

\title{The Brownian Lévy transformation is exact}
\author{Zhi Ling\thanks{National University of Singapore. Email: \texttt{lingzhi@nus.edu.sg}.}}
\date{}

\usepackage{xurl}

\begin{document}
\maketitle

\begin{abstract}
We prove that the Lévy transformation of Brownian motion is exact. This settles the long-standing question of whether the transformation is ergodic. The proof combines an interval cutting argument with an estimate for the first entrance of the iterated paths into a small neighborhood of zero at an independent exponential time.
\end{abstract}

\section{Introduction}

Let $\mathcal{W}:=C_0([0,\infty),\mathbb{R})$, let $\mathbb{P}$ be Wiener measure, and let $\mathcal{F}$ be the completion of the Borel $\sigma$-field of $\mathcal{W}$. The coordinate process $B_t(\omega):=\omega(t)$ is a standard Brownian motion. Let $L$ denote the local time of $B$ at zero in Tanaka's normalization. The Lévy transformation is the map $T:\mathcal{W}\to\mathcal{W}$ given by
\begin{equation}
(TB)_t:=\int_0^t\operatorname{sgn}(B_s)\,\mathrm{d}B_s=|B_t|-L_t.
\label{eq.levy_transform}
\end{equation}
Since the stochastic integral in \eqref{eq.levy_transform} has quadratic variation $t$, Lévy's characterization implies that $T$ preserves $\mathbb{P}$. Our main result considers the information retained by its iterations.

\begin{theorem}
\label{thm.exactness}
The Brownian Lévy transformation is exact. Equivalently,
$$
\bigcap_{n\geq0}T^{-n}\mathcal{F}
$$
is trivial (modulo $\mathbb{P}$). Therefore, $T$ is strongly mixing and ergodic.
\end{theorem}

\subsection{Background and related work}

The Lévy transformation originates in Lévy's 1940 work \cite{Levy1940}. Its preservation of Wiener measure naturally raises the question of its behavior under repeated iteration. The ergodicity problem, raised by Marc Yor \cite{DubinsSmorodinsky1992}, was already known in the late 1970s \cite{DubinsEmeryYor1993} and was recorded as an open problem in the first edition of Revuz and Yor's monograph \cite[p.~257]{RevuzYor1991}.

Discrete analogues illustrate how ergodicity depends on the choice of discrete transformation. Dubins and Smorodinsky proved that a modified discrete Lévy transformation is Bernoulli \cite{DubinsSmorodinsky1992}, whereas Fujita exhibited a nonergodic version \cite{Fujita2008}. Collevecchio, Hamza, Shi, and Williams later gave a necessary and sufficient condition for ergodicity of general adapted transformations of simple symmetric random walks \cite{CollevecchioHamzaShiWilliams2022}.

In the Brownian setting, Dubins, Émery, and Yor linked ergodicity to the characterization of continuous martingales invariant in law under the Lévy transformation \cite{DubinsEmeryYor1993}. Zero sets provide a concrete description of this dependence. Malric proved that $T^n B$ and $T^m B$ almost surely have common zeros at positive times when $m\geq n+2$ \cite{Malric1995}. Donati-Martin, Shi, and Yor studied the joint law of the last zeros before time one of Brownian motion and its Lévy transform, finding a strictly positive covariance \cite{DonatiMartinShiYor2000}.

Research moved from proving that iterated paths are dense to developing criteria for mixing. Malric proved that the union of the zero sets of the iterates is almost surely dense in time \cite{Malric2003}, and later that almost every orbit is dense in Wiener space for uniform convergence on compact intervals \cite{Malric2012}. Brossard and Leuridan gave another proof \cite{BrossardLeuridan2012}. Prokaj established sufficient conditions for strong mixing in terms of how rapidly the values of the iterated paths at time one approach zero \cite{Prokaj2013}. His subsequent work strengthened this approach by providing a sufficient condition for exactness \cite{Prokaj2014}.

Section~\ref{sec.reduction} reduces the theorem to the entrance estimate according to Prokaj's criterion. Section~\ref{sec.rooted_records} develops the change of root and the record estimates. Section~\ref{sec.entrance_estimate} derives logarithmic integrability and proves the entrance estimate by excursion compensation.

\subsection{Notation}

Throughout the paper we use $B^{(n)}:=T^nB$ and
$$
Z_n(t):=\min_{0\leq j<n}|B_t^{(j)}|.
$$
For $j\geq0$, let $L^{(j)}$ be the local time of $B^{(j)}$ at zero, so that
\begin{equation}
B_t^{(j+1)}=|B_t^{(j)}|-L_t^{(j)}.
\label{eq.iterate_tanaka}
\end{equation}
A superscript in parentheses is an iteration index, and a time subscript is Brownian time. The value of $\operatorname{sgn}(0)$ is immaterial. Write $\mathbb{E}:=\mathbb{E}_{\mathbb{P}}$, $x\wedge y:=\min\{x,y\}$, $x\vee y:=\max\{x,y\}$, and $x_+:=\max\{x,0\}$. We use $\inf\varnothing:=\infty$ and empty products equal to one. In estimates, $C$ is independent of the displayed parameters and may change between occurrences. The notation $w_{\cdot\wedge t}$ means the entire path $s\mapsto w_{\min\{s,t\}}$, held constant after $t$.

\section{A reduction to entrance locations}
\label{sec.reduction}

We use the following consequence of \cite[Theorems~8 and~9]{Prokaj2013} and \cite[Theorem~1]{Prokaj2014}.

\begin{proposition}[Prokaj's criterion]\label{prop.prokaj_criterion}
\begin{equation}
\liminf_{n\to\infty}\frac{Z_{n+1}(1)}{Z_n(1)}=0\quad\text{almost surely} \quad \Longrightarrow \quad T \text{~is exact.}
\label{eq.ratio_condition}
\end{equation}
Since $\liminf_n Z_{n+1}(1)/Z_n(1) \in \{0,1\}$ a.s., we only need to show the liminf is smaller than 1 a.s.
\end{proposition}

We will first work at a random observation time. Let $\tau$ be an exponential random variable of rate $1/2$, independent of the Brownian path. We continue to write $\mathbb{P}$ for the probability on the enlarged space and $\mathbb{E}$ for its expectation. Malric's orbit density theorem \cite{Malric2012} implies that $Z_n(t)\downarrow0$ for every $t\geq0$ almost surely. So for every $r>0$ we may define the finite entrance index and its entrance value by
$$
\nu(r):=\inf\{n\geq0:|B_\tau^{(n)}|\leq r\},\qquad
U(r):=|B_\tau^{(\nu(r))}|.
$$
Thus $\nu(r)$ counts iterations, and $U(r)$ is the value at the first entrance into $[0,r]$.

\begin{lemma}[Entrance estimate]
\label{lem.entrance_estimate}
There is a constant $C$ such that
\begin{equation}
\mathbb{P}\bigl(  U(r) \in [(1-\varepsilon)r, r]\bigr)
\leq C\varepsilon\log(e/r),\qquad 0<r\leq1,\quad 0<\varepsilon\leq\tfrac12.
\label{eq.entrance_bound}
\end{equation}
\end{lemma}

We will prove Lemma~\ref{lem.entrance_estimate} in Section~\ref{sec.entrance_estimate}. Here we show it forces proportional decreases of $Z_n(\tau)$. The first step is to condition the endpoint on the local times and the interval left by the preceding zeros. We first describe the information revealed by the local times and the zeros of the preceding iterates.

\subsection{The endpoint in a random interval}

Set $X:=B_\tau$ and $\lambda_j:=L_\tau^{(j)}$. With these local time values fixed, define functions of a real variable $x$ by
$$
f_0(x):=x,\qquad f_{j+1}(x):=|f_j(x)|-\lambda_j,\qquad j\geq0.
$$
Tanaka's identity \eqref{eq.iterate_tanaka} gives $f_j(X)=B_\tau^{(j)}$. Thus all the iterated values at time $\tau$ can be recovered from the single endpoint $X$ and the local time values. Each $B_\tau^{(j)}$ has a continuous distribution. We work on the event of probability one on which none of these values is zero.

For $n\geq1$, remove from $\mathbb{R}$ all zeros of $f_0,\ldots,f_{n-1}$, and let $I_n$ be the remaining open interval that contains $X$. Define
$$
\mathcal{G}_n:=\sigma(\lambda_0,\ldots,\lambda_{n-1},I_n),\qquad
\mathcal{G}_\infty:=\sigma\left(\bigcup_{n\geq1}\mathcal{G}_n\right).
$$
A random interval is specified by its two endpoints, which may be infinite. The $\sigma$-fields $\mathcal{G}_n$ are increasing. Indeed, $\mathcal{G}_{n+1}$ gives the local times needed to construct $f_0,\ldots,f_{n-1}$, and $I_n$ is the unique component of their common complement that contains $I_{n+1}$.

On $I_n$, none of $f_0,\ldots,f_{n-1}$ changes sign. Consequently, $f_n$ is affine there, with slope $1$ or $-1$. Given $\mathcal{G}_n$, we know whether $f_n$ has a zero in $I_n$ and, if so, where that zero lies. Such a zero cuts $I_n$ into two intervals. The next interval $I_{n+1}$ is the one that contains $X$. If there is no zero in $I_n$, then $I_{n+1}=I_n$.

\begin{lemma}
\label{lem.interval_density}
For every $n\geq1$,
\begin{equation}
Z_n(\tau)=\operatorname{dist}(X,\partial I_n).
\label{eq.interval_distance}
\end{equation}
Conditionally on $\mathcal{G}_n$, if $I_n=(\alpha,\alpha+\ell)$ is bounded, then $X-\alpha$ has density
\begin{equation}
\rho_{\ell,p}(u):=\frac{p\sinh(\ell-u)+(1-p)\sinh u}{\cosh\ell-1},\qquad 0<u<\ell,
\label{eq.interval_density}
\end{equation}
where $p\in[0,1]$ is $\mathcal{G}_n$-measurable. If $I_n$ is a half-line, the distance of $X$ from its finite endpoint has density $e^{-u}$ on $(0,\infty)$. Moreover, $|I_n|\downarrow0$ almost surely, and $X$ is $\mathcal{G}_\infty$-measurable.
\end{lemma}

\begin{proof}
We first prove the distance identity. It holds for $n=1$, since the only zero removed is that of $f_0(x)=x$. Suppose it holds at index $n$, and write $z:=\operatorname{dist}(X,\partial I_n)$. The function $f_n$ is affine with slope of absolute value one on $I_n$. If $|f_n(X)|<z$, its affine zero lies in $I_n$ at distance $|f_n(X)|$ from $X$. If $|f_n(X)|\geq z$, the $1$-Lipschitz property of $f_n$ prevents any of its zeros from being closer than $z$. In either case,
$$
\operatorname{dist}(X,\partial I_{n+1})
=\min\{z,|f_n(X)|\}.
$$
Induction and $f_j(X)=B_\tau^{(j)}$ give \eqref{eq.interval_distance}.

We next determine the conditional density. Write $\lambda:=(\lambda_0,\ldots,\lambda_{n-1})$, and let $\mu$ denote the joint law of $(\lambda,X)$. In this argument, the zero sets $f_j(x)=0$ are viewed as subsets of the joint $(\lambda,x)$ space. Choose a smooth compactly supported function $\varphi(\lambda,x)$ whose support avoids all these sets for $j<n$. Apply Itô's formula to
$$
e^{-t/2}\varphi(L_t^{(0)},\ldots,L_t^{(n-1)},B_t).
$$
The term involving $\mathrm{d}L_t^{(j)}$ vanishes because this measure is supported on $B_t^{(j)}=0$, where $\varphi$ and its derivatives vanish. The initial value also vanishes, since $B_0=0$. After taking expectations and letting the terminal time tend to infinity, we obtain
$$
\begin{aligned}
0
&=\mathbb{E}\int_0^\infty e^{-t/2}
\left(\tfrac12\partial_{xx}\varphi-\tfrac12\varphi\right)
(L_t^{(0)},\ldots,L_t^{(n-1)},B_t)\,\mathrm{d}t\\
&=\int(\partial_{xx}\varphi-\varphi)\,\mathrm{d}\mu.
\end{aligned}
$$
The second equality uses the density $\frac12e^{-t/2}$ of $\tau$.

Let $\mu_\lambda$ be the conditional law of $X$ given $\lambda$. In a product neighborhood that avoids all the zero sets, take $\varphi(\lambda,x):=\psi(\lambda)\chi(x)$, with both factors smooth and compactly supported. Varying $\psi$ gives
$$
\int\bigl(\chi''(x)-\chi(x)\bigr)\,\mu_\lambda(\mathrm{d}x)=0
$$
for almost every $\lambda$ in that neighborhood. Use a countable cover by such neighborhoods and a countable dense family of test functions in each one. The resulting identities then hold simultaneously on every interval left after the zeros are removed. On each such interval, the conditional measure satisfies $\mu_\lambda''=\mu_\lambda$ in the sense of distributions. Every distributional solution of this equation is a linear combination of $e^x$ and $e^{-x}$. Thus the conditional measure has a smooth density satisfying the same differential equation.

For fixed $\lambda$, the functions $f_0,\ldots,f_{n-1}$ have only finitely many zeros. Conditioning further on $I_n$ therefore selects one of the finitely many intervals and normalizes the conditional density on it. After translation to $(0,\ell)$, every nonnegative solution of the density equation has the form
$$
a\sinh(\ell-u)+b\sinh u,\qquad a,b\geq0.
$$
The coefficients are nonnegative because the density has nonnegative limits at both endpoints. Each basis function has integral $\cosh\ell-1$. Normalization gives \eqref{eq.interval_density}, with $p:=a/(a+b)$. On a half-line, integrability removes the increasing exponential. Normalization leaves $e^{-u}$ for the distance from the finite endpoint.

Finally, we show that the intervals shrink to $X$. Fix $c>0$. If $\ell\geq c$, then
$$
\rho_{\ell,p}(u)
\leq\frac{\sinh\ell}{\cosh\ell-1}
=\coth(\ell/2)
\leq\coth(c/2).
$$
The density on a half-line obeys the same upper bound. The part of an interval within distance $v$ of its boundary has total length at most $2v$. By \eqref{eq.interval_distance},
$$
\mathbb{P}\bigl(|I_n|\geq c,\ Z_n(\tau)\leq v\bigr)
\leq2v\coth(c/2),\qquad v>0.
$$
The intervals are nested, whereas $Z_n(\tau)\downarrow0$. Thus
$$
\begin{aligned}
\mathbb{P}\bigl(\lim_n|I_n|\geq c\bigr)
&\leq\limsup_n\mathbb{P}\bigl(Z_n(\tau)>v\bigr)+2v\coth(c/2)\\
&=2v\coth(c/2).
\end{aligned}
$$
Letting $v\downarrow0$, and then taking all positive rational $c$, proves $|I_n|\downarrow0$. Since every $I_n$ contains $X$, its finite endpoints converge to $X$. These endpoints are $\mathcal{G}_\infty$-measurable, so $X$ is also $\mathcal{G}_\infty$-measurable.
\end{proof}

\subsection{From entrance locations to proportional decreases}

For the rest of this subsection, write
$$
z_n:=Z_n(\tau),\qquad \delta_n:=1-\frac{z_{n+1}}{z_n},\qquad n\geq1.
$$
The variables $z_n$ are positive almost surely. By \eqref{eq.interval_distance}, $\delta_n$ is the proportion by which the distance from $X$ to the interval boundary decreases when $I_n$ is replaced by $I_{n+1}$.

If $I_n$ is bounded and $f_n$ has a zero inside it, let $J_n$ be the shorter of the two intervals produced by the cut. Choose the left interval if their lengths are equal. Set $J_n:=\varnothing$ when there is no cut or $I_n$ is unbounded. Given $\mathcal{G}_n$, the interval $J_n$ is known, but we do not yet know whether it contains $X$. Define
$$
\pi_n:=\mathbb{P}(X\in J_n\mid\mathcal{G}_n).
$$
The event $\{X\in J_n\}$ belongs to $\mathcal{G}_{n+1}$, because the latter reveals which interval was selected.

\begin{lemma}
\label{lem.interval_cutting}
If Lemma~\ref{lem.entrance_estimate} holds, then $\delta_n\geq1/2$ infinitely often, almost surely.
\end{lemma}

\begin{proof}
Suppose, to the contrary, that $A:=\{\delta_n<1/2\text{ eventually}\}$ has positive probability. We will first show that $\sum_n\pi_n<\infty$ on $A$. We will then use this summability to bound the decreases of $z_n$, and show that this bound contradicts the entrance estimate.

Fix an integer $N\geq1$ and set $A_N:=\{\delta_n<1/2\text{ for all }n\geq N\}$. Lemma~\ref{lem.interval_density} shows that $X$, and hence every $\delta_n$, is $\mathcal{G}_\infty$-measurable. Martingale convergence therefore gives
$$
\mathbb{P}(A_N\mid\mathcal{G}_{n+1})\longrightarrow\mathbf{1}_{A_N}
\qquad\text{almost surely}.
$$
We compare this limit with what happens when the shorter interval is selected. Choose coordinates so that $I_n=(0,\ell)$ and its shorter component lies on the left. Write this component as $J_n=(0,a)$, where $a\leq\ell/2$, and continue to denote the endpoint in these coordinates by $X$. On the event $\{X\in J_n\}$, we have $I_{n+1}=(0,a)$. If $X=u\in(3a/4,7a/8)$, then its distance from the old boundary is $u$, and its distance from the new boundary is $a-u$. Therefore
$$
\delta_n=1-\frac{a-u}{u}=2-\frac au>\frac23.
$$
If $a\leq1$, the conditional probability that $X$ lies in $(3a/4,7a/8)$ given $\mathcal{G}_{n+1}$ is at least $1/64$. To see this, apply \eqref{eq.interval_density} to the selected interval $(0,a)$. On the indicated subinterval, both $\sinh(a-u)$ and $\sinh u$ are at least $a/8$, while $\cosh a-1\leq a^2$. The conditional density is therefore at least $1/(8a)$ on an interval of length $a/8$.

Since $|I_n|\to0$, we have $a\leq1$ at every sufficiently late cut. Whenever the shorter interval is selected at such an index $n\geq N$, the preceding bound gives
$$
\mathbb{P}(A_N\mid\mathcal{G}_{n+1})\leq1-\frac1{64}.
$$
On $A_N$, the conditional probabilities on the left converge to one. The shorter interval can therefore be selected only finitely often on $A_N$. Since $A=\bigcup_N A_N$, the same conclusion holds on $A$. The conditional Borel-Cantelli lemma, applied to $\{X\in J_n\}\in\mathcal{G}_{n+1}$, now gives $\sum_n\pi_n<\infty$ on $A$.

Next we bound the probability of a decrease in terms of $\pi_n$. Keep the coordinates $I_n=(0,\ell)$ and $J_n=(0,a)$, and assume $\ell\leq1$. Let $\delta_n(x)$ denote the relative decrease when $X=x$. For $x>a$,
$$
\min\{x-a,\ell-x\}\geq(1-a/x)\min\{x,\ell-x\}.
$$
The two minima are the new and old boundary distances, respectively, so $\delta_n(x)\leq a/x$. The same bound holds for $0<x<a$, since $\delta_n(x)\leq1<a/x$. Consequently,
$$
\{x:\delta_n(x)\geq u\}\subset(0,\min\{a/u,\ell\}],
\qquad 0<u\leq1.
$$

Let $\mu_n$ be the conditional law of the endpoint in these coordinates given $\mathcal{G}_n$. The density \eqref{eq.interval_density} implies
$$
\mu_n((0,s])\leq C(s/a)^2\mu_n((0,a]),\qquad a\leq s\leq\ell.
$$
Indeed, for $0\leq v\leq1$,
$$
v\leq\sinh v\leq(\sinh1)v.
$$
It is therefore enough to compare the integrals of $\ell-x$ and $x$ over $(0,s]$ and $(0,a]$. The corresponding ratios are at most $2s/a$ and $(s/a)^2$. Since $s/a\geq1$, both are bounded by $2(s/a)^2$, and the bound is preserved under a nonnegative mixture. Taking $s:=\min\{a/u,\ell\}$ and using $\mu_n((0,a])=\pi_n$ gives
$$
\mathbb{P}(\delta_n\geq u\mid\mathcal{G}_n)\leq C\pi_n/u^2.
$$
This bound also holds when there is no cut, because then $\delta_n=0$.

We now use this tail bound on a fixed part of $A$ with positive probability. By a countable union argument, there exist an event $A'\subset A$ with $\mathbb{P}(A')>0$ and deterministic constants $n_0\geq1$, $K_0<\infty$, and $0<r_0\leq1$ such that, on $A'$,
$$
z_{n_0}\geq r_0,\qquad
\sum_{n\geq n_0}\pi_n\leq K_0,\qquad
|I_n|\leq1\ \text{and}\ \delta_n<1/2\quad\text{for every }n\geq n_0.
$$
Such constants exist because $z_n>0$, the intervals shrink, and the sum of the $\pi_n$ is finite on $A$.

Fix $0<\theta\leq1/2$. Integrating the conditional tail bound gives
$$
\begin{aligned}
\mathbb{E}\left[\delta_n\mathbf{1}_{\{\delta_n\geq\theta\}}\mid\mathcal{G}_n\right]
&=\theta\,\mathbb{P}(\delta_n\geq\theta\mid\mathcal{G}_n)
+\int_\theta^1\mathbb{P}(\delta_n\geq u\mid\mathcal{G}_n)\,\mathrm{d}u\\
&\leq C\pi_n/\theta
\end{aligned}
$$
on $\{|I_n|\leq1\}$. The conditional expectation is also at most one. Hence
$$
\theta\,\mathbb{E}\left[\delta_n\mathbf{1}_{\{\delta_n\geq\theta\}}\mid\mathcal{G}_n\right]
\leq\min\{\theta,C\pi_n\}.
$$
We cannot condition directly on $A'$, since it may depend on the entire sequence. Instead, for $n\geq n_0$, define the $\mathcal{G}_n$-measurable event
$$
E_n:=\left\{|I_n|\leq1,\ \sum_{k=n_0}^n\pi_k\leq K_0\right\}.
$$
Since $A'\subset E_n$, conditional expectation yields
$$
\begin{aligned}
\theta\,\mathbb{E}[\mathbf{1}_{A'}\delta_n\mathbf{1}_{\{\delta_n\geq\theta\}}]
&\leq\mathbb{E}\left[\mathbf{1}_{E_n}\theta\,
\mathbb{E}[\delta_n\mathbf{1}_{\{\delta_n\geq\theta\}}\mid\mathcal{G}_n]\right]\\
&\leq\mathbb{E}[\mathbf{1}_{E_n}\min\{\theta,C\pi_n\}].
\end{aligned}
$$
Summing over $n$ gives
\begin{equation}
\theta\,\mathbb{E}\left[\mathbf{1}_{A'}\sum_{n\geq n_0}\delta_n\mathbf{1}_{\{\delta_n\geq\theta\}}\right]
\leq\mathbb{E}\sum_{n\geq n_0}\mathbf{1}_{E_n}\min\{\theta,C\pi_n\}
\longrightarrow0
\label{eq.large_cut_bound}
\end{equation}
as $\theta\downarrow0$. For this last limit, each summand tends to zero, and
$$
\sum_{n\geq n_0}\mathbf{1}_{E_n}\pi_n\leq K_0.
$$
Dominated convergence first applies to the sum for each sample and then to its expectation.

The final step compares \eqref{eq.large_cut_bound} with the entrance estimate over many thresholds $r$. If $z_{n+1}<z_n$ and $r\in[z_{n+1},z_n)$, the first iterate with absolute value at most $r$ is the one that produces $z_{n+1}$. Thus
$$
U(r)=z_{n+1},\qquad 1-\frac{U(r)}r\leq\delta_n.
$$
In particular, $U(r)<(1-\theta)r$ requires $\delta_n\geq\theta$. On $A'$, the interval $[z_{n+1},z_n)$ has logarithmic length
$$
\int_{z_{n+1}}^{z_n}\frac{\mathrm{d}r}{r}
=-\log(1-\delta_n)\leq2\delta_n,\qquad n\geq n_0.
$$
These intervals cover $(0,r_0)$ apart from their endpoints, since $z_{n_0}\geq r_0$ and $z_n\downarrow0$. It follows from \eqref{eq.large_cut_bound} that
$$
\theta\,\mathbb{E}\left[\mathbf{1}_{A'}\int_0^{r_0}
\mathbf{1}_{\{U(r)<(1-\theta)r\}}\,\frac{\mathrm{d}r}{r}\right]
\leq2\theta\,\mathbb{E}\left[\mathbf{1}_{A'}\sum_{n\geq n_0}
\delta_n\mathbf{1}_{\{\delta_n\geq\theta\}}\right]
\longrightarrow0.
$$

The entrance estimate gives a positive lower bound for the same quantity. Write $\alpha:=\mathbb{P}(A')$, let $C$ be the constant in \eqref{eq.entrance_bound}, and set
$$
r_\theta:=\exp(1-\alpha/(2C\theta)).
$$
For sufficiently small $\theta$, we have $r_\theta<r_0$. For every $r\in[r_\theta,r_0]$, Lemma~\ref{lem.entrance_estimate} gives
$$
\mathbb{P}(A',\ U(r)<(1-\theta)r)
\geq\alpha-C\theta\log(e/r)\geq\frac\alpha2.
$$
Tonelli's theorem therefore implies
$$
\theta\,\mathbb{E}\left[\mathbf{1}_{A'}\int_0^{r_0}
\mathbf{1}_{\{U(r)<(1-\theta)r\}}\,\frac{\mathrm{d}r}{r}\right]
\geq\frac{\alpha\theta}{2}\log\frac{r_0}{r_\theta}
\longrightarrow\frac{\alpha^2}{4C}>0.
$$
This contradicts the preceding limit and proves the lemma.
\end{proof}

Assuming Lemma~\ref{lem.entrance_estimate}, we can now prove Theorem~\ref{thm.exactness}. Lemma~\ref{lem.interval_cutting} gives
$$
\liminf_n\frac{Z_{n+1}(\tau)}{Z_n(\tau)}\leq\frac12
\qquad\text{almost surely}.
$$
Brownian scaling commutes with $T$, and the common scaling factor cancels in each ratio. The probability of the corresponding event at a deterministic time $t>0$ is therefore the same for every $t$. Since $\tau$ is independent of the paths, integrating this probability against $\frac12e^{-t/2}\,\mathrm{d}t$ gives one. The constant probability must consequently be one. In particular, the liminf at time one is at most $1/2$ almost surely. Proposition~\ref{prop.prokaj_criterion} gives exactness.

\section{Records at zeros of the iterated paths}
\label{sec.rooted_records}

The entrance estimate requires us to examine several iterates during one excursion. On an interval where a path stays away from zero, its sign and local time are constant. Tanaka's formula then expresses the next iterate as an affine function of that path. We use these affine relations to identify the heights at which other iterates can reach zero.

This section is arranged in three stages. We first define a probability measure for observing the system at a zero of a specified iterate. We then record the heights at which higher iterates can reach zero. Finally, we move the reference zero from one iterate to another. This change of reference will give estimates for the number of records and the number of iterations between them.

\subsection{The probability measure at a zero}

The proof sometimes requires examining iterates with indices smaller than that of the reference path. We therefore work on the natural extension of $T$, which provides a stationary sequence of iterates indexed by all integers. Let $\Omega:=\mathcal{W}^{\mathbb{Z}}$, with its product Borel $\sigma$-field. To specify the law of a finite block $(B^{(i)},\ldots,B^{(j)})$, give $B^{(i)}$ Wiener measure and obtain each successive coordinate by applying $T$. These laws are consistent because $T$ preserves Wiener measure. They define a probability measure on $\Omega$ that is invariant under a common shift of the iteration indices and satisfies
$$
TB^{(j)}=B^{(j+1)},\qquad j\in\mathbb{Z}.
$$
We again denote this measure by $\mathbb{P}$ and write $\mathcal{A}$ for the completed product $\sigma$-field. The nonnegative coordinates have exactly their original joint law. For each $j\in\mathbb{Z}$, let $L^{(j)}$ be the local time at zero of $B^{(j)}$ in Tanaka's normalization. Then \eqref{eq.iterate_tanaka} holds at every integer index. We may also adjoin the independent variable $\tau$.

Let
$$
\mathcal{H}_t:=\sigma(B_s^{(j)}:j\in\mathbb{Z},\ 0\leq s\leq t),
$$
with the usual augmentation. Every coordinate is Brownian with respect to this common filtration. To verify this, fix a coordinate and a finite collection of coordinates whose histories we wish to observe. Choose the lowest index among them. Every other coordinate in the collection is adapted to the Brownian filtration of that lowest coordinate and is a stochastic integral with quadratic variation $t$. Lévy's characterization makes it Brownian in that filtration. Its future increments are therefore independent of the observed finite collection of past histories. A monotone class argument extends this independence to $\mathcal{H}_t$. The strong Markov property and excursion compensation can consequently be used while keeping the past of all lower iterates in the conditioning history.

We will repeatedly use the following affine relation. Suppose $i<j$ and $B^{(i+1)},\ldots,B^{(j-1)}$ stay nonzero on a time interval. Their signs and local times are constant there. Repeated application of Tanaka's formula gives, for two times $u,s$ in that interval,
\begin{equation}
B_s^{(j)}-B_u^{(j)}
=\left(\prod_{k=i+1}^{j-1}\operatorname{sgn}(B_u^{(k)})\right)
\bigl(B_s^{(i+1)}-B_u^{(i+1)}\bigr).
\label{eq.affine_propagation}
\end{equation}
We also use this relation to prove that the intervening paths stay nonzero. Starting at $u$, compute their affine expressions under the assumption that no sign changes. If none of these expressions reaches zero on the proposed interval, a first zero of an intervening path would contradict \eqref{eq.affine_propagation} up to that time. The same reasoning can be applied from a terminal time by considering the last possible zero before it.

We also use the reflection identities
$$
L_t^{(i)}=-\inf_{0\leq s\leq t}B_s^{(i+1)},\qquad
|B_t^{(i)}|=B_t^{(i+1)}-\inf_{0\leq s\leq t}B_s^{(i+1)}.
$$
For the first identity, Tanaka's formula gives $B_s^{(i+1)}\geq-L_s^{(i)}\geq-L_t^{(i)}$ for $s\leq t$. On the support of $\mathrm{d}L^{(i)}$, it gives $B^{(i+1)}=-L^{(i)}$. Since local time is continuous and increases from zero to $L_t^{(i)}$, the running infimum must also be at most $-L_t^{(i)}$. The second identity then follows from Tanaka's formula. Thus local time is the absolute value of the driving iterate's running minimum.

We now choose observation times using local time. On $\Omega\times[0,\infty)$, with $\sigma$-field $\mathcal{A}\otimes\mathcal{B}([0,\infty))$, define
\begin{equation}
\mathbb{Q}(A):=\mathbb{E}\int_0^\infty e^{-t/2}\mathbf{1}_A(\omega,t)\,\mathrm{d}L_t^{(0)}.
\label{eq.zero_measure}
\end{equation}
Local time increases only when $B^{(0)}$ is zero, so $\mathbb{Q}$ observes the system at such a zero. The exponential factor makes its total mass finite and, with the chosen rate and normalization, equal to one:
$$
\mathbb{Q}(\Omega\times[0,\infty))
=\int_0^\infty\frac{e^{-t/2}}{\sqrt{2\pi t}}\,\mathrm{d}t=1.
$$
Here we used $\mathbb{E}L_t^{(0)}=\sqrt{2t/\pi}$. Thus $\mathbb{Q}$ is a probability measure. A sample consists of the paths and an observation time $t$ with $B_t^{(0)}=0$. We call the pair $(0,t)$ the root and write $\mathbb{E}_{\mathbb{Q}}$ for expectation under this measure. The first entry in a root is an iteration index, and the second is a Brownian time.

For an integer $p$ and a time $g$, define the history
$$
H_{p,g}:=\left(\bigl((B_{s\wedge g}^{(p+j)})_{s\geq0}\bigr)_{j\in\mathbb{Z}},\ g\right).
$$
This history contains every path up to time $g$, holds the paths constant after $g$, and relabels original coordinate $p+j$ as coordinate $j$. If $B_g^{(p)}=0$, its reference path is therefore zero at its observation time. We shift only the iteration indices. Brownian time still starts at zero. For any function $F$ of such a history, write $F_{p,g}:=F(H_{p,g})$, keeping any other arguments of $F$ unchanged. Index stationarity gives, for every nonnegative measurable $F$,
$$
\mathbb{E}\int_0^\infty e^{-g/2}F(H_{p,g})\,\mathrm{d}L_g^{(p)}
=\mathbb{E}_{\mathbb{Q}}F(H_{0,t}).
$$
This identity allows the same measure $\mathbb{Q}$ to describe observations at zeros of any chosen iterate. On the space of stopped histories, we also write $\mathbb{Q}$ for the law of $H_{0,t}$ induced by \eqref{eq.zero_measure}.

\begin{lemma}
\label{lem.disjoint_local_times}
For distinct integers $i,j$, the measure $\mathrm{d}L^{(i)}$ gives zero mass to the zero set of $B^{(j)}$. Consequently, the local time measures of distinct coordinates are mutually singular.
\end{lemma}

\begin{proof}
We prove that each local time measure gives zero mass to the other path's zero set, by induction on $j-i>0$. First restrict attention to an open interval on which the intermediate paths $B^{(i+1)},\ldots,B^{(j-1)}$ are nonzero. By \eqref{eq.affine_propagation} and Tanaka's formula, on this interval
$$
B^{(j)}=v|B^{(i)}|-vL^{(i)}+c,
$$
where $v\in\{-1,1\}$ and $c$ are constant. At a common zero of $B^{(i)}$ and $B^{(j)}$, this identity requires $L^{(i)}=c/v$. The measure $\mathrm{d}L^{(i)}$ gives zero mass to the set where the continuous increasing function $L^{(i)}$ has any one fixed value. Thus it gives zero mass to these common zeros.

For the reverse direction, the times at which Brownian local time has a fixed value form an interval, possibly a single point. The support of Brownian local time is the Brownian zero set, so the interior of a nontrivial such interval contains no zeros. There are therefore at most two zeros at any fixed local time level. The common zeros under consideration consequently form a set of at most two points. Since $L^{(j)}$ is continuous, $\mathrm{d}L^{(j)}$ also gives them zero mass.

There are only countably many open components on which all intermediate paths are nonzero. Outside these components, at least one intermediate path is zero. By the induction hypothesis, this exceptional set has zero mass for both $\mathrm{d}L^{(i)}$ and $\mathrm{d}L^{(j)}$. This proves the induction step. When $j=i+1$, there are no intermediate paths, and the same affine argument proves the initial case. Finally, each local time measure is supported on the zero set of its own path. The claim just proved therefore implies mutual singularity.
\end{proof}

\subsection{Zero levels and their records}

Work under $\mathbb{Q}$ at the root $(0,t)$. We want to know how the higher iterates respond when the absolute value of the reference path moves away from zero. Hold all local time values at $t$ fixed and use Tanaka's formula to express the higher iterates as functions of this absolute value. Near zero, each expression is affine. It remains valid as long as no earlier expression changes sign.

Suppose one of these affine expressions has value $x\neq0$ at the root and slope $\sigma\in\{-1,1\}$. Its zero is at the signed level $b:=-\sigma x$. A positive level lies in the direction in which the reference absolute value increases. A negative level means that this affine expression moves away from zero in that direction. Applying one more Lévy transform with local time value $\lambda$ changes the zero level to
$$
-\sigma\operatorname{sgn}(x)(|x|-\lambda)
=b-\operatorname{sgn}(b)\lambda.
$$
For the first iterate the expression is $u-L_t^{(0)}$, where $u$ is the reference absolute value. Thus the signed zero levels are given recursively by
\begin{equation}
b_1:=L_t^{(0)},\qquad
b_{j+1}:=b_j-\operatorname{sgn}(b_j)L_t^{(j)},\quad j\geq1.
\label{eq.zero_level_recursion}
\end{equation}
In particular, $|b_j|=|B_t^{(j)}|$. Lemma~\ref{lem.disjoint_local_times} implies that all these levels are nonzero, $\mathbb{Q}$-almost surely.

For $a>0$, define the forward waiting index by
$$
\kappa(a):=\inf\{j\geq1:0<b_j\leq a\}.
$$
This is the first higher iteration index whose positive zero level is at most $a$. In a history $H_{p,g}$, the same definition gives
$$
\kappa_{p,g}(a):=\kappa(a)(H_{p,g}).
$$
The resulting index is measured relative to $p$, so the corresponding original coordinate is $p+\kappa_{p,g}(a)$.

We will use the positive levels that are smaller than every earlier positive level. Define their set by
$$
\mathcal{R}:=\left\{b_j:j\geq1,\quad
0<b_j<\inf_{\substack{1\leq k<j\\b_k>0}}b_k\right\}.
$$
We call these levels records. A record at height $y$ has all earlier positive zero levels above $y$. The corresponding affine expressions therefore remain nonzero while the reference absolute value increases from zero to $y$. This is why records identify zeros that can be reached before an earlier iterate changes sign.

For $y\in\mathcal{R}$, define its iteration index and its gap to the next smaller record by
$$
p_y:=\inf\{j\geq1:b_j=y\},\qquad
D(y):=\inf\{j>p_y:0<b_j<y\}-p_y.
$$
For a Borel set $I\subset(0,\infty)$, let $N(I):=\#(\mathcal{R}\cap I)$. Record heights decrease as their iteration indices increase. Thus a record higher than $y$ occurs before $y$.

The first record is $b_1=L_t^{(0)}>0$, $\mathbb{Q}$-almost surely. Starting from its index, which is one, add the gaps after every record above $a$. This gives the index of the first record at or below $a$:
\begin{equation}
\kappa(a)=1+\sum_{\substack{y\in\mathcal{R}\\y>a}}D(y).
\label{eq.record_decomposition}
\end{equation}
For now this is an identity in $[1,\infty]$. If no positive level at most $a$ occurs, there are either infinitely many records above $a$ or a final record above $a$ with infinite gap. The right-hand side is infinite in either case. We will prove finiteness after establishing the change of root formula.

\subsection{Changing the reference zero}

Fix a record $y\in\mathcal{R}$ at the source root $(0,t)$, and put $p:=p_y$. At time $t$, coordinate $p$ has absolute value $y$. Let its last zero before $t$, its sign at $t$, and its absolute excursion up to $t$ be
$$
g_y:=\sup\{0\leq s<t:B_s^{(p)}=0\},\qquad
\eta_y:=\operatorname{sgn}(B_t^{(p)}),\qquad
e_y(s):=|B_{(g_y+s)\wedge t}^{(p)}|,\quad s\geq0.
$$
The path $e_y$ starts at zero and is held at $y$ after elapsed time $t-g_y$. We will replace the source root $(0,t)$ by the earlier zero $(p,g_y)$. The proof below will show that $t$ is the first time this excursion reaches absolute height $y$.

To identify the source from the new root, we need to inspect lower iterates. Consider a history rooted at $(0,g)$ and define
$$
\beta_i:=-|B_g^{(-i)}|\prod_{k=1}^{i-1}\operatorname{sgn}(B_g^{(-k)}),
\qquad i\geq1.
$$
While coordinates $-1,\ldots,-i$ remain nonzero, their local times stay constant. Solving their affine relations in terms of the reference path gives
$$
|B_s^{(-i)}|=|\beta_i|-\operatorname{sgn}(\beta_i)B_s^{(0)}.
$$
Thus $\beta_i$ is the signed value of the reference path at which this affine expression reaches zero. During a reference excursion of sign $\eta\in\{-1,1\}$, a lower coordinate can reach zero by absolute height $a$ only if $0<\eta\beta_i\leq a$. Define
$$
\kappa_{\mathrm b}(a,\eta):=\inf\{i\geq1:0<\eta\beta_i\leq a\},\qquad
\overline{\kappa}_{\mathrm b}(a):=\frac12\sum_{\eta\in\{-1,1\}}\kappa_{\mathrm b}(a,\eta).
$$
The search visits coordinates $-1,-2,\ldots$ in that order. The bar denotes an average over the two excursion signs, with the history fixed. Both quantities are nonincreasing in $a$. Set $\kappa_{\mathrm b}(0,\eta):=\infty$. A subscript $p,g$ means that this search is performed in $H_{p,g}$, so relative coordinate $-i$ is original coordinate $p-i$.

We next describe the excursion measure used in the change of root formula. Let $\mathbf n_{\mathrm{pos}}$ be Itô's measure of positive Brownian excursions in Tanaka's local time normalization \cite[Chapter~XII]{RevuzYor1999}. Excursion compensation expresses an expected sum over excursions as an integral of their contributions against $\mathbf n_{\mathrm{pos}}$ and the local time at their starting zeros. In this sense, $\mathbf n_{\mathrm{pos}}$ gives the rate of positive excursions per unit of local time. It is a measure rather than a probability law. Negative excursions have the same measure for their absolute paths.

Write $\zeta(e)$ for an excursion's lifetime and
$$
T_a(w):=\inf\{s\geq0:w(s)=a\}
$$
for a path's first hitting time of $a$. We omit the path argument when it is clear. Let $\mathbb{E}_x$ denote expectation for Brownian motion started at $x$. The discounted hitting probability on $(0,y)$ solves $u''=u$ with $u(0)=0$ and $u(y)=1$. Hence
$$
\mathbb{E}_x[e^{-T_y/2}\mathbf{1}_{\{T_y<T_0\}}]
=\frac{\sinh x}{\sinh y},\qquad 0<x<y.
$$
In the chosen normalization, passage from Brownian motion started at $x$ to positive excursion measure divides this expression by $2x$ and lets $x\downarrow0$. We obtain
$$
\mathbf n_{\mathrm{pos}}(e^{-T_y/2}\mathbf{1}_{\{T_y<\zeta\}})
=\frac1{2\sinh y}.
$$
It follows that
$$
\Pi_y(F):=2\sinh y\,\mathbf n_{\mathrm{pos}}
\bigl(e^{-T_y/2}F(e_{\cdot\wedge T_y})\mathbf{1}_{\{T_y<\zeta\}}\bigr)
$$
defines a probability law on positive excursion paths stopped on first reaching $y$. Here $\mathbf n_{\mathrm{pos}}(F)$ means integration against excursion measure. The stopped path is held at $y$ after $T_y$, and $F$ may use its entire history before that time.

The following identity describes the expected sum over records after moving to their new roots. On its right-hand side, the new history has law $\mathbb{Q}$, the sign is equally likely to be $1$ or $-1$, and the absolute excursion has law $\Pi_y$. The height is integrated with density $1/\sinh y$. This separation of the history, sign, and absolute excursion is the main reason for changing the reference zero.

\begin{proposition}[Change of root]
\label{prop.marked_rerooting}
For every nonnegative measurable function $\Phi$ of a height, a stopped history, a sign, and a stopped positive path,
\begin{equation}
\begin{aligned}
&\mathbb{E}_{\mathbb{Q}}\sum_{y\in\mathcal{R}}\Phi(y,H_{p_y,g_y},\eta_y,e_y)\\
&\qquad=\int_0^\infty\frac1{2\sinh y}\sum_{\eta\in\{-1,1\}}
\mathbb{E}_{\mathbb{Q}}\int\Phi(y,H_{0,t},\eta,e)\,\Pi_y(\mathrm{d}e)\,\mathrm{d}y.
\end{aligned}
\label{eq.marked_rerooting}
\end{equation}
\end{proposition}

\begin{proof}
The main step is a correspondence between two descriptions of the same event. In the first description, we observe a zero of a lower iterate and a record among its higher iterates. In the second, we observe an excursion of the record iterate when it first reaches the record height. We prove that these descriptions agree and that source local time becomes ordinary measure on the height variable. Excursion compensation will then give the stated formula.

Throughout this correspondence, fix $W:=B^{(0)}$ as the path that will provide the record. Its possible source paths have indices $-1,-2,\ldots$. We first record two almost sure properties that will be used for all its excursions. The union of the zero sets of these lower paths is dense in time. Indeed, for any fixed rational open interval, stationarity makes the probability of having no zero of coordinates $-n,\ldots,-1$ equal to the probability of having no zero of coordinates $0,\ldots,n-1$. The latter tends to zero by Malric's theorem on the density of the union of the zero sets \cite{Malric2003}. Taking all rational intervals proves the claim.

Also, every nontrivial excursion of $W$ starts at a time when all lower coordinates are nonzero. For a fixed lower coordinate and excursions reaching height $1/k$, excursion compensation reduces the expected number of exceptions on a bounded time interval to an integral against $\mathrm{d}L^{(0)}$ over the zero set of that lower coordinate. Lemma~\ref{lem.disjoint_local_times} makes this integral zero. Taking all lower indices, positive integers $k$, and bounded time intervals covers every nontrivial excursion.

Consider such an excursion of $W$, starting at time $g$, with sign $\eta$, absolute path $e$, and lifetime $\zeta$. Write
$$
\max e:=\max_{0\leq s\leq\zeta}e(s).
$$
Fix $0<y\leq\max e$, and set
$$
t_y:=g+T_y(e),\qquad m:=\kappa_{\mathrm b,0,g}(y,\eta).
$$
The signed levels $\beta_i$ in the following argument are computed in the history $H_{0,g}$. We claim that the source associated with height $y$ is $(-m,t_y)$.

First, $m$ is finite. While the paths with indices $-1,\ldots,-i+1$ remain nonzero, affine propagation gives the displacement
$$
B_{g+s}^{(-i+1)}-B_g^{(-i+1)}
=-\eta\operatorname{sgn}(\beta_i)e(s).
$$
Until coordinate $-i$ itself reaches zero, its absolute value is therefore
$$
|\beta_i|-\eta\operatorname{sgn}(\beta_i)e(s).
$$
This expression can vanish only if the excursion reaches the positive height $\eta\beta_i$. If $m$ were infinite, no such height would lie in $(0,y]$. The affine formula would then keep coordinate $-1$ nonzero until $t_y$. Applying the same argument successively to coordinates $-2,-3,\ldots$ would keep every lower path nonzero on $(g,t_y)$. This contradicts the density of their zeros. The same argument shows that all coordinates $-i$ with $i<m$ remain nonzero on $[g,t_y]$.

By the definition of $m$, we have $0<\eta\beta_m\leq y$, so $\operatorname{sgn}(\beta_m)=\eta$. The displacement of the driving path $B^{(-m+1)}$ is consequently $-e(s)$. Let
$$
M_s:=\sup_{0\leq v\leq s}e(v).
$$
Tanaka's reflection identity says that $L^{(-m)}$ increases by the amount by which this driving path goes below its previous running minimum. Its initial distance from that minimum is $|B_g^{(-m)}|=|\beta_m|$. Hence
\begin{equation}
L_{g+s}^{(-m)}-L_g^{(-m)}=(M_s-|\beta_m|)_+,
\qquad 0\leq s\leq T_y(e).
\label{eq.source_local_time}
\end{equation}
At $s=T_y(e)$, both $e(s)$ and $M_s$ equal $y$. The driving path is at its running minimum, so $B_{t_y}^{(-m)}=0$. Thus $(-m,t_y)$ is a source zero.

We next check that $W$ is a record relative to this source. Its value at $t_y$ is $\eta y$, and its affine slope as a function of the source absolute value is $-\eta$. Its zero level is therefore $y$. For an intervening coordinate $-i$, where $i<m$, its absolute value at $t_y$ is
$$
|\beta_i|-\eta\operatorname{sgn}(\beta_i)y,
$$
and its slope as a function of the source absolute value is $\eta\operatorname{sgn}(\beta_i)$. Its zero level is consequently
$$
y-\eta\beta_i.
$$
Since $i<m$, either $\eta\beta_i<0$ or $\eta\beta_i>y$. In the first case this level exceeds $y$, and in the second it is negative. Every preceding positive level is therefore larger than $y$. This proves that $W$ supplies a record at relative index $m$ from the source $(-m,t_y)$.

We now prove the converse. Fix $m\geq1$, and consider $\mathrm{d}L^{(-m)}$-almost every time $t$ at which $W$ supplies a record of height $y$ from $(-m,t)$. Define
$$
v:=\prod_{k=-m+1}^{-1}\operatorname{sgn}(B_t^{(k)}),\qquad
d_u:=B_u^{(-m+1)}-B_t^{(-m+1)},\quad 0\leq u\leq t.
$$
At a source zero, $B_t^{(-m+1)}=-L_t^{(-m)}$ is the running minimum of the driving path. Therefore $d_u\geq0$ and $d_t=0$. Also,
$$
d_0=L_t^{(-m)}\geq y,
$$
because $L_t^{(-m)}$ is the first positive zero level from this source and $y$ is a record. Let $g<t$ be the last time at which $d_g=y$. Such a time exists by continuity. For $g<u\leq t$, we have $0\leq d_u<y$.

Every positive zero level before the record exceeds $y$. Applying the affine relation backwards from $t$ therefore keeps the intervening coordinates nonzero on $[g,t]$ and gives
$$
W_u-W_t=vd_u,\qquad W_t=-vy,\qquad W_g=0.
$$
It follows that $g$ is the last zero of $W$ before $t$ and
$$
|W_u|=y-d_u,\qquad g\leq u\leq t.
$$
Thus the excursion containing $t$ has absolute value at most $y$ before $t$.

For almost every $t$ under consideration, it first reaches $y$ at $t$. Indeed, if $|W_u|=y$ for some $g<u<t$, then $d_u=0$. The driving path has already reached its terminal running minimum at $u$, so $L^{(-m)}$ is constant on $[u,t]$. A continuous increasing function has only countably many levels with a nontrivial interval of constancy, and its Stieltjes measure gives zero mass to each such level. These exceptional times therefore have zero $\mathrm{d}L^{(-m)}$-mass. Writing $e$ for the absolute excursion from $g$, we obtain $t=g+T_y(e)$.

Let $\eta$ be the sign of this excursion, and compute the $\beta_i$ in $H_{0,g}$. The intervening coordinates remain nonzero until the first visit to height $y$, so $\eta\beta_i\notin(0,y]$ for $i<m$. For coordinate $-m$, the driver decreases by $y$ between $g$ and $t$ and ends at its running minimum. Hence its initial distance from the running minimum is at most $y$, and the displacement is towards that minimum. In terms of the signed level, this says
$$
0<\eta\beta_m=|\beta_m|\leq y.
$$
We have recovered $m=\kappa_{\mathrm b,0,g}(y,\eta)$, as required.

It remains to identify the measures on the two sides of this correspondence. Fix an excursion starting at $g$ and a possible source index $-m$. If $\eta\beta_m>0$, the heights assigned to this source index are exactly
$$
\eta\beta_m\leq y<
\inf_{\substack{i<m\\\eta\beta_i>0}}\eta\beta_i,
$$
intersected with $(0,\max e]$. If $\eta\beta_m<0$, no height is assigned to it. On the interior of the indicated interval, \eqref{eq.source_local_time} identifies source local time increments with increments of $M$. For a continuous running maximum, the change of variables is
$$
\int_0^\zeta F_{g+s}\,\mathrm{d}M_s
=\int_0^{\max e}F_{g+T_y(e)}\,\mathrm{d}y
$$
for every nonnegative measurable function of time $F$, with the paths fixed. Each increase of the maximum corresponds to the first visit to a new height. The heights at which $M$ has a nontrivial plateau form a countable set, and plateaus carry no $\mathrm{d}M$-mass, so they do not affect this identity.

Restrict this change of variables to the heights assigned to $-m$. The endpoints of the height intervals and the excursion maxima can also be omitted. Over all excursions and indices, their first hitting times form a countable set, which has zero mass for every continuous local time. Let $C_m(t)$ denote the event that $W$ supplies a record from $(-m,t)$. The correspondence and the change of variables now give
\begin{equation}
\int_0^\infty\mathbf{1}_{C_m(t)}F_t\,\mathrm{d}L_t^{(-m)}
=\sum_g\int_0^{\max e_g}
\mathbf{1}_{\{\kappa_{\mathrm b,0,g}(y,\eta_g)=m\}}
F_{g+T_y(e_g)}\,\mathrm{d}y,
\label{eq.record_height_measure}
\end{equation}
where $g$ runs over the starting times of excursions of $W$, and $e_g,\eta_g$ are their absolute paths and signs.

We can now apply excursion compensation. For a time $t$ inside an excursion of $W$, let $g$ be its starting time and define
$$
\Phi_t:=\Phi\left(|W_t|,H_{0,g},\operatorname{sgn}(W_t),
\bigl(|W_{(g+s)\wedge t}|\bigr)_{s\geq0}\right).
$$
In the sum on the left of \eqref{eq.marked_rerooting}, first group records by their iteration index $m$. Shift all indices by $-m$, so that the record path becomes $W=B^{(0)}$ and the source path becomes $B^{(-m)}$. Index stationarity then gives the first equality below. The second uses \eqref{eq.record_height_measure}, and the third is excursion compensation:
$$
\begin{aligned}
&\mathbb{E}_{\mathbb{Q}}\sum_{y\in\mathcal{R}}\Phi(y,H_{p_y,g_y},\eta_y,e_y)\\
&\quad=\sum_{m\geq1}\mathbb{E}\int_0^\infty
e^{-t/2}\mathbf{1}_{C_m(t)}\Phi_t\,\mathrm{d}L_t^{(-m)}\\
&\quad=\mathbb{E}\sum_g e^{-g/2}\int_0^{\max e_g}
e^{-T_y(e_g)/2}
\Phi\bigl(y,H_{0,g},\eta_g,(e_g)_{\cdot\wedge T_y(e_g)}\bigr)\,\mathrm{d}y\\
&\quad=\mathbb{E}\int_0^\infty e^{-t/2}\left[
\sum_{\eta\in\{-1,1\}}\int_0^\infty
\mathbf n_{\mathrm{pos}}\bigl(e^{-T_y/2}
\Phi(y,H_{0,t},\eta,e_{\cdot\wedge T_y})
\mathbf{1}_{\{T_y<\zeta\}}\bigr)\,\mathrm{d}y
\right]\mathrm{d}L_t^{(0)}.
\end{aligned}
$$
In the second equality, summing over $m$ removes the indicator because every reached height has exactly one source index. The discount splits as $e^{-t/2}=e^{-g/2}e^{-T_y/2}$. In the last line, the history at the excursion's starting zero is kept fixed while its sign and path are integrated against excursion measure. Substituting the definitions of $\Pi_y$ and $\mathbb{Q}$ proves \eqref{eq.marked_rerooting}.

The same argument also proves finiteness of the backward waiting indices under $\mathbb{Q}$. Fix a rational $a>0$ and a sign $\eta$. We have shown that every excursion of that sign reaching $a$ starts from a history in which $\kappa_{\mathrm b}(a,\eta)<\infty$. Applying compensation to the exceptional histories gives
$$
\frac{\mathbb{Q}(\kappa_{\mathrm b}(a,\eta)=\infty)}{2\sinh a}=0.
$$
Taking all positive rational $a$ and both signs, then using monotonicity in $a$, proves finiteness for every $a>0$ on a single event of full $\mathbb{Q}$-measure.
\end{proof}

We will also need almost sure properties to remain valid after changing the root. Let $E$ be a set of stopped histories with $\mathbb{Q}$-measure zero. Apply \eqref{eq.marked_rerooting} to
$$
\Phi(y,H,\eta,e):=\mathbf{1}_{[a,b]}(y)\mathbf{1}_E(H).
$$
For $0<a<b<\infty$, this gives
$$
\mathbb{E}_{\mathbb{Q}}\sum_{y\in\mathcal{R}\cap[a,b]}
\mathbf{1}_E(H_{p_y,g_y})=0.
$$
Covering $(0,\infty)$ by countably many such intervals shows that, almost surely, no record root has a history in $E$. Thus any fixed almost sure statement about the reference history can be used simultaneously at all record roots.

\begin{lemma}[Record indices]
For almost every rooted sample, every record $y$, with $p:=p_y$, satisfies
\begin{equation}
p=\kappa_{\mathrm b,p,g_y}(y,\eta_y).
\label{eq.source_index}
\end{equation}
For each fixed $a>0$, a record $y\leq a$ is the first positive level at most $a$ if and only if
\begin{equation}
\kappa_{\mathrm b,p,g_y}(y,\eta_y)
<\kappa_{\mathrm b,p,g_y}(a-y,-\eta_y).
\label{eq.first_record_filter}
\end{equation}
Its gap to the next smaller record satisfies
\begin{equation}
D(y)=\kappa_{p,g_y}(y)
\label{eq.record_gap}
\end{equation}
as an identity in $[1,\infty]$.
\end{lemma}

\begin{proof}
The correspondence in Proposition~\ref{prop.marked_rerooting} identifies the source index by searching backwards from the record root. This gives \eqref{eq.source_index}.

For the second assertion, compute the backward levels $\beta_i$ in $H_{p,g_y}$. The proof of that correspondence shows that the zero levels preceding $y$ in the source history are
$$
y-\eta_y\beta_i,\qquad i<p.
$$
Since $y$ is already a record, each of these levels is either negative or larger than $y$. Such a level belongs to $(0,a]$ precisely when
$$
0<-\eta_y\beta_i\leq a-y.
$$
The search with sign $-\eta_y$ and threshold $a-y$ detects the first index with this property. It cannot stop at $i=p$, because the source index satisfies $\eta_y\beta_p>0$. Thus there is no positive level at most $a$ before $y$ exactly when this opposite-sign search stops after $p$. Substituting \eqref{eq.source_index} gives \eqref{eq.first_record_filter}. If $a=y$, the opposite-sign threshold is zero and its waiting index is infinite, as stipulated in the definition.

To prove the gap identity, write $g:=g_y$ and define
$$
b_i':=b_i(H_{p,g}),\qquad d:=\kappa_{p,g}(y).
$$
The unprimed levels refer to the source history at $(0,t)$, and the primed levels refer to the new history at $(p,g)$. We first exclude equalities between $y$ and the primed levels. For any fixed new history, the set $\{b_i':i\geq1\}$ is countable. Since the change of root formula integrates $y$ against a density, it gives zero expected count to record summands for which $y=b_i'$ for some $i$. We may therefore assume that no such equality occurs. The primed levels are also nonzero by Lemma~\ref{lem.disjoint_local_times} and the transfer of almost sure properties to record roots.

By the definition of $d$, we have $b_i'\notin[0,y]$ for $i<d$. During the excursion from $g$ to $t$, the absolute value of coordinate $p$ stays in $[0,y]$ and first reaches $y$ at $t$. Their zero levels lie outside the range of the reference absolute value, so the affine relations keep coordinates $p+i$, $1\leq i<d$, nonzero throughout this interval. Their local times, and the local time of coordinate $p$, remain constant.

We can now compare the two recursions. Since $b_p=y$,
$$
b_{p+1}=y-L_g^{(p)}=y-b_1'.
$$
For $i<d$, the numbers $b_i'$ and $y-b_i'$ have opposite signs, because $b_i'<0$ or $b_i'>y$. Applying \eqref{eq.zero_level_recursion} to both histories, with the same local time values, gives inductively
$$
b_{p+i}=y-b_i',\qquad 1\leq i\leq d.
$$
If $d=\infty$, this assertion is understood for every finite $i$. For $i<d$, the source level $y-b_i'$ lies outside $(0,y)$. If $d<\infty$, then $0<b_d'<y$, so $y-b_d'\in(0,y)$. The next smaller source record is therefore at index $p+d$. If $d=\infty$, no smaller source record occurs. This proves \eqref{eq.record_gap} in both cases.
\end{proof}

\begin{corollary}[Record intensity and mean waiting index]
For Borel $I\subset(0,\infty)$ and nonnegative measurable $F$ on $(0,\infty)\times[1,\infty]$,
\begin{equation}
\mathbb{E}_{\mathbb{Q}}N(I)=\int_I\frac{\mathrm{d}y}{\sinh y},
\label{eq.record_intensity}
\end{equation}
and
\begin{equation}
\mathbb{E}_{\mathbb{Q}}\sum_{y\in\mathcal{R}}F(y,D(y))
=\int_0^\infty\frac{\mathbb{E}_{\mathbb{Q}}F(y,\kappa(y))}{\sinh y}\,\mathrm{d}y.
\label{eq.waiting_rerooting}
\end{equation}
Moreover,
\begin{equation}
m(a):=\mathbb{E}_{\mathbb{Q}}\kappa(a)=\coth(a/2),\qquad a>0.
\label{eq.mean_waiting_bound}
\end{equation}
\end{corollary}

\begin{proof}
In \eqref{eq.marked_rerooting}, take $\Phi$ to be the indicator that the height belongs to $I$. Each $\Pi_y$ has total mass one, and summing the two signs gives \eqref{eq.record_intensity}. Next use \eqref{eq.record_gap} to write the gap at a source record as the forward waiting index in its new history. The same change of root formula gives \eqref{eq.waiting_rerooting}.

We prove the mean formula by truncating the waiting index first. For integers $n\geq1$, set
$$
m_n(a):=\mathbb{E}_{\mathbb{Q}}[\kappa(a)\wedge n].
$$
The decomposition \eqref{eq.record_decomposition} implies
$$
\kappa(a)\wedge n\leq1+\sum_{\substack{y\in\mathcal{R}\\y>a}}(D(y)\wedge n).
$$
Taking expectations and applying \eqref{eq.waiting_rerooting}, we obtain
$$
m_n(a)\leq1+\int_a^\infty\frac{m_n(y)}{\sinh y}\,\mathrm{d}y.
$$
Substitute this inequality into its own right-hand side repeatedly. The $k$-fold integral is over ordered heights $a<y_1<\cdots<y_k$, so it is $1/k!$ times the corresponding integral over all $k$ heights. Since $m_n\leq n$, the remaining term after $k$ substitutions is bounded by
$$
\frac{n}{k!}\left(\int_a^\infty\frac{\mathrm{d}y}{\sinh y}\right)^k,
$$
which tends to zero. The terms generated by the constant one sum to an exponential series. Thus
$$
m_n(a)\leq\exp\left(\int_a^\infty\frac{\mathrm{d}y}{\sinh y}\right)
=\coth(a/2).
$$
Monotone convergence as $n\to\infty$ gives $m(a)\leq\coth(a/2)<\infty$.

We may now take expectations in the untruncated decomposition:
$$
m(a)=1+\int_a^\infty\frac{m(y)}{\sinh y}\,\mathrm{d}y.
$$
Repeated substitution in this equality, followed by omission of the nonnegative remaining term, bounds $m(a)$ below by every partial sum of the same exponential series. Hence $m(a)\geq\coth(a/2)$, proving \eqref{eq.mean_waiting_bound}.

For each rational $a>0$, finiteness of the mean gives $\kappa(a)<\infty$ almost surely. Monotonicity then makes all $\kappa(a)$ finite on one event of full measure. On this event, record heights decrease to zero and only finitely many records lie above any fixed positive threshold.
\end{proof}

Before estimating record counts, we specify which information is unaffected by changing an excursion sign. Fix a record height $y$, the new history $H_{p,g_y}$, and the absolute excursion $e_y$. Write $g:=g_y$. During this excursion, the local time of coordinate $p$ is constant, so
$$
B_{g+s}^{(p+1)}=e_y(s)-L_g^{(p)},\qquad 0\leq s\leq t-g.
$$
This expression does not use the excursion sign $\eta_y$. It determines the path of coordinate $p+1$ up to $t$. Applying Tanaka's formula successively then determines every higher path and its local time up to $t$, still without using $\eta_y$.

Starting from the source level $b_p=y$, the recursion \eqref{eq.zero_level_recursion} now determines all subsequent source levels and their indices relative to $p$. It also determines the last zeros of the later record coordinates and the histories of higher coordinates at those zeros. In particular, the later record heights and gaps do not depend on $\eta_y$ when the new history and absolute excursion are fixed. The original index $p$ can change with the sign, but the relative indices and gaps do not.

This observation combines with a simple consequence of \eqref{eq.marked_rerooting}. Let $G(y,H,e)\geq0$ be a weight that does not depend on the excursion sign. Suppose an event $A(y,H,\eta)$ holds for at least one of the two signs for each height and history under consideration. Put
$$
G_y:=G(y,H_{p_y,g_y},e_y),\qquad
A_y:=A(y,H_{p_y,g_y},\eta_y).
$$
On the right-hand side of the change of root formula, the signs have equal weights. Since at least one sign satisfies $A$, their average obeys
$$
\frac12\sum_{\eta\in\{-1,1\}}G(y,H,e)\mathbf{1}_{A(y,H,\eta)}
\geq\frac12G(y,H,e).
$$
Applying the formula to both sides yields
\begin{equation}
\mathbb{E}_{\mathbb{Q}}\sum_{y\in\mathcal{R}}G_y\mathbf{1}_{A_y}
\geq\frac12\mathbb{E}_{\mathbb{Q}}\sum_{y\in\mathcal{R}}G_y.
\label{eq.sign_selection}
\end{equation}
We will use this inequality to select the first record in an interval while retaining at least half of the expected weight.

\subsection{Record positions and counts}

Fix an interval $I:=(a,2a]$. For a record $y\in I$, define
$$
r_I(y):=N(I\cap(y,\infty)).
$$
This is the number of records in $I$ that precede $y$. Thus $r_I(y)=0$ means that $y$ is the first record in $I$.

At the new root of $y$, the two waiting indices $\kappa_{\mathrm b}(y,1)$ and $\kappa_{\mathrm b}(y,-1)$ are finite and distinct. They are distinct because the same signed level $\beta_i$ cannot qualify for both signs. Choose $\eta$ for which $\kappa_{\mathrm b}(y,\eta)$ is smaller. Since $y>a$, we have $2a-y\leq y$, and monotonicity gives
$$
\kappa_{\mathrm b}(y,\eta)
<\kappa_{\mathrm b}(y,-\eta)
\leq\kappa_{\mathrm b}(2a-y,-\eta).
$$
The criterion for the first record \eqref{eq.first_record_filter}, with threshold $2a$, shows that this sign makes $y$ the first record in $I$. Consequently, \eqref{eq.sign_selection} applies to the event $\{r_I(y)=0\}$ for every weight on records in $I$ that does not depend on the excursion sign.

We will turn this selection probability into bounds on counts. The elementary fact we need is the following. Let $V$ be an integrable nonnegative integer-valued variable, and suppose that, for some $0<\theta<1$,
$$
\mathbb{P}(V>m)\geq\theta\,\mathbb{E}(V-m)_+,\qquad m\in\{0,1,2,\ldots\}.
$$
Set $s_m:=\mathbb{E}(V-m)_+$. Since $V$ is integer-valued,
$$
s_{m+1}=s_m-\mathbb{P}(V>m)\leq(1-\theta)s_m.
$$
Iteration gives
\begin{equation}
\mathbb{E}(V-m)_+\leq(1-\theta)^m\mathbb{E}V.
\label{eq.selection_tail_bounds}
\end{equation}

\begin{lemma}[Single records]
For $I:=(a,2a]$,
\begin{equation}
\mathbb{E}_{\mathbb{Q}}\sum_{y\in\mathcal{R}\cap I}(r_I(y)+1)D(y)
\leq2\mathbb{E}_{\mathbb{Q}}\sum_{y\in\mathcal{R}\cap I}D(y),
\label{eq.record_position_bound}
\end{equation}
and, for every integer $m\geq0$,
\begin{equation}
\mathbb{E}_{\mathbb{Q}}(N(I)-m)_+\leq(\log2)2^{-m}.
\label{eq.single_count_tail}
\end{equation}
\end{lemma}

\begin{proof}
For a record $y\in I$, let
$$
S_y:=\kappa(a)-p_y
=\sum_{\substack{z\in\mathcal{R}\cap I\\z\leq y}}D(z).
$$
This is the total number of iterations from $y$ until the first positive level at most $a$. The second expression shows that it depends only on the gap after $y$ and the later gaps in $I$. By the reconstruction above, it does not depend on the sign of the excursion attached at the new root of $y$.

If we sum $S_y$ over all records in $I$, a particular gap $D(z)$ is included once for $z$ itself and once for every preceding record in $I$. Hence
$$
\sum_{y\in\mathcal{R}\cap I}S_y
=\sum_{z\in\mathcal{R}\cap I}(r_I(z)+1)D(z).
$$
At the first record in $I$, the weight $S_y$ includes every gap in $I$. Sign selection therefore gives
$$
\begin{aligned}
\mathbb{E}_{\mathbb{Q}}\sum_{y\in\mathcal{R}\cap I}S_y
&\leq2\mathbb{E}_{\mathbb{Q}}\sum_{y\in\mathcal{R}\cap I}
S_y\mathbf{1}_{\{r_I(y)=0\}}\\
&=2\mathbb{E}_{\mathbb{Q}}\sum_{y\in\mathcal{R}\cap I}D(y).
\end{aligned}
$$
This proves \eqref{eq.record_position_bound}. The last expectation is finite, since \eqref{eq.record_decomposition} identifies it with $m(a)-m(2a)$.

For the count bound, fix $m\geq0$ and give a record weight one if at least $m$ records in $I$ follow it, and weight zero otherwise. This weight again depends only on later records, so it does not depend on the attached excursion sign. If there are $N(I)$ records, exactly $(N(I)-m)_+$ of them have weight one. At the first record the weight is one precisely when $N(I)>m$. Applying sign selection gives
$$
\mathbb{Q}(N(I)>m)\geq\frac12\mathbb{E}_{\mathbb{Q}}(N(I)-m)_+.
$$
Finally, the intensity formula gives
$$
\mathbb{E}_{\mathbb{Q}}N(I)
=\int_a^{2a}\frac{\mathrm{d}y}{\sinh y}
\leq\int_a^{2a}\frac{\mathrm{d}y}{y}
=\log2.
$$
Apply \eqref{eq.selection_tail_bounds} with $\theta=1/2$ to obtain \eqref{eq.single_count_tail}.
\end{proof}

We also need to count records seen after two successive changes of root. For Borel sets $I,J\subset(0,\infty)$ bounded away from zero, define
$$
N(I,J):=\sum_{y\in\mathcal{R}\cap I}
N_{p_y,g_y}\bigl(J\cap(y,\infty)\bigr).
$$
To form a counted pair, first choose a record $y\in I$ in the source history. Then move to its root $(p_y,g_y)$ and choose a record $h\in J$ of that new history with $h>y$. We call $y$ the parent and $h$ the child. The restriction $h>y$ ensures that the child coordinate stays away from zero during the parent excursion, whose absolute height is at most $y$.

\begin{lemma}[Pairs of records]
For $I:=(a,2a]$, $J:=(b,2b]$, and integers $m\geq0$,
\begin{equation}
\mathbb{E}_{\mathbb{Q}}(N(I,J)-m)_+
\leq(\log2)^2(3/4)^m.
\label{eq.paired_record_tail}
\end{equation}
\end{lemma}

\begin{proof}
We will apply the same counting argument as for single records. This time, we select the first pair by selecting its parent first and then its child. Each sign selection retains at least half of the expected weight. To justify applying the argument twice, we must show that the weight placed on a pair does not depend on either of the two excursion signs.

First, the mean count is finite. Change the root at the parent and apply the record intensity formula in the resulting history:
$$
\begin{aligned}
\mathbb{E}_{\mathbb{Q}}N(I,J)
&=\int_I\frac{\mathrm{d}y}{\sinh y}
\int_{J\cap(y,\infty)}\frac{\mathrm{d}h}{\sinh h}\\
&\leq(\log2)^2.
\end{aligned}
$$

Next we specify the order of the pairs. For a counted pair $(y,h)$, write $p:=p_y$ and $g:=g_y$, and let $q$ be the original iteration index of its child. Since $h$ is a record above $y$ in the parent history, it occurs before the first positive level at most $y$ in that history. Thus, by \eqref{eq.record_gap},
$$
p<q<p+\kappa_{p,g}(y)=p+D(y).
$$
The last expression is the index of the next parent record in the source history. All children of $y$ therefore have indices between $p$ and the next parent index. Children of different parents occupy disjoint intervals of iteration indices, in the same order as their parents. We order all counted pairs by their child indices $q$.

We now determine what is needed to recover the pairs following a fixed pair $(y,h)$. Let
$$
\gamma:=\sup\{0\leq s<g:B_s^{(q)}=0\},\qquad
e_h(s):=|B_{(\gamma+s)\wedge g}^{(q)}|,\quad s\geq0.
$$
Thus $e_h$ is the child's absolute excursion from its last zero before $g$ up to its first visit to $h$ at time $g$. Fix the child history $H_{q,\gamma}$ and the two absolute excursions $e_h$ and $e_y$. Their durations determine
$$
g=\gamma+T_h(e_h),\qquad t=g+T_y(e_y).
$$
During the parent excursion, all coordinates between $p$ and $q$ stay nonzero. The child has positive zero level $h$ in the parent history, so its absolute value decreases by the parent absolute excursion value. Hence
\begin{equation}
|B_s^{(q)}|=
\begin{cases}
e_h(s-\gamma),&\gamma\leq s\leq g,\\
h-e_y(s-g),&g\leq s\leq t.
\end{cases}
\label{eq.parent_child_path}
\end{equation}
Since $h>y$, this path is positive throughout $(\gamma,t]$. Its local time is constant there. The history up to $\gamma$ and the absolute path in \eqref{eq.parent_child_path} therefore determine $B^{(q+1)}$ up to $t$ by Tanaka's formula. Repeating the formula determines all higher paths and their local times. None of this uses the sign of the child excursion or the sign of the parent excursion.

There are two kinds of pairs that can follow $(y,h)$. First come the remaining children of the same parent. At time $g$, the child's zero level in the parent history is $h$. Starting from this level, \eqref{eq.zero_level_recursion} and the now known higher local times recover every later level in that history. We keep the subsequent records in $J$ whose heights exceed $y$. This recovers all remaining children of the current parent.

The other pairs belong to later parents. The comparison of levels in the proof of \eqref{eq.record_gap} gives the source level at coordinate $q$ as
$$
b_q=y-h<0.
$$
At time $t$, start the recursion from this level, with current positive record $y$. It recovers all later source records and identifies the later parents in $I$. Their coordinates exceed $q$, so their paths up to $t$ are already known. We can therefore find their last zeros and all histories of higher coordinates needed to compute their children. This recovers every pair after $(y,h)$, including their relative order.

Consequently, with $y$, $h$, $H_{q,\gamma}$, $e_h$, and $e_y$ fixed, the list of later pairs does not depend on either excursion sign. The original iteration labels may change, but the relative order and the number of later pairs remain fixed. In particular, for each integer $m\geq0$, the weight
$$
G_m(y,h):=\mathbf{1}_{\{\text{at least }m\text{ counted pairs follow }(y,h)\}}
$$
has this property.

We now perform the two sign selections. In the remaining sums, $(y,h)$ always ranges over the pairs counted by $N(I,J)$. Define
$$
A_y:=\{N(I\cap(y,\infty))=0\},\qquad
A_h:=\{N_{p,g}(J\cap(h,\infty))=0\}.
$$
The event $A_y$ means that $y$ is the first parent record in $I$. The event $A_h$ means that $h$ is the first record in $J$ in the parent history. Because $h>y$, every record preceding $h$ in that history also exceeds $y$. Thus, for a counted pair, $A_h$ is exactly the condition that its child is the first counted child of that parent.

Fix $m$ and write $G:=G_m(y,h)$ in each summand. Apply the change of root formula at the parent. With its new history and absolute excursion fixed, $A_h$ is determined by that history, and $G$ does not depend on the parent excursion sign. At least one parent sign satisfies $A_y$, by the selection of the first record in $I$. Hence
$$
\mathbb{E}_{\mathbb{Q}}\sum_{(y,h)}G\mathbf{1}_{A_y}\mathbf{1}_{A_h}
\geq\frac12\mathbb{E}_{\mathbb{Q}}\sum_{(y,h)}G\mathbf{1}_{A_h}.
$$

For the second selection, keep the parent height $y$ and absolute excursion $e_y$ fixed. Denote the weight reconstructed from a child root by $\widetilde G_{y,e_y}(h,H,e_h)$. The reconstruction above shows that this weight does not depend on the child excursion sign. Averaging the parent sign in \eqref{eq.marked_rerooting} gives
$$
\begin{aligned}
&\mathbb{E}_{\mathbb{Q}}\sum_{(y,h)}G\mathbf{1}_{A_h}\\
&\quad=\int_I\frac{\mathrm{d}y}{\sinh y}\int\Pi_y(\mathrm{d}e_y)\,
\mathbb{E}_{\mathbb{Q}}\sum_{h\in\mathcal{R}\cap J\cap(y,\infty)}
\widetilde G_{y,e_y}(h,H_{p_h,g_h},e_h)
\mathbf{1}_{\{N(J\cap(h,\infty))=0\}}.
\end{aligned}
$$
In the inner expectation, the parent history is the reference history with law $\mathbb{Q}$. The parent excursion $e_y$ is an independent parameter with law $\Pi_y$. We may therefore apply the change of root formula again, now to its child records. For $h\in J=(b,2b]$, the inequality $2b-h\leq h$ ensures that at least one child sign satisfies the criterion for the first record in $J$. Sign selection gives
$$
\mathbb{E}_{\mathbb{Q}}\sum_{(y,h)}G\mathbf{1}_{A_h}
\geq\frac12\mathbb{E}_{\mathbb{Q}}\sum_{(y,h)}G.
$$

There is at most one pair satisfying both $A_y$ and $A_h$. If present, it is the first counted pair, because its parent is first and its child is first within that parent. Its weight $G_m$ is at most $\mathbf{1}_{\{N(I,J)>m\}}$. On the other hand, the sum of $G_m$ over all pairs is $(N(I,J)-m)_+$. Combining the two selections yields
$$
\mathbb{Q}(N(I,J)>m)
\geq\frac14\mathbb{E}_{\mathbb{Q}}(N(I,J)-m)_+.
$$
The counting bound \eqref{eq.selection_tail_bounds}, with $\theta=1/4$, and the mean estimate at the start of the proof give \eqref{eq.paired_record_tail}.
\end{proof}

\section{Proof of the entrance estimate}
\label{sec.entrance_estimate}

We now return to the first entrance at the independent exponential time $\tau$. The eventual excursion calculation will involve a product of two quantities from the same root history. A backward waiting index bounds how many reference iterates can contribute, and a record count bounds how many higher iterates can produce an entrance. We need to control the expectation of this product.

The count estimates from the previous section give exponential tails. We first bound $\mathbb{E}_{\mathbb{Q}}[\kappa(a)\log\kappa(a)]$ for $a>0$. The change of root formula will then express the needed product in terms of a forward waiting index and record counts, to which these two estimates can be applied.

\subsection{Logarithmic integrability of the waiting index}

\begin{lemma}
With $m(a)=\coth(a/2)$ as in \eqref{eq.mean_waiting_bound},
\begin{equation}
\mathbb{E}_{\mathbb{Q}}[\kappa(a)\log\kappa(a)]
\leq C m(a)(1+\log m(a)),\qquad a>0.
\label{eq.waiting_entropy}
\end{equation}
\end{lemma}

\begin{proof}
For integers $n\geq1$, define
$$
h_n(a):=\mathbb{E}_{\mathbb{Q}}
[(\kappa(a)\wedge n)\log(\kappa(a)\wedge n)].
$$
We will find a bound that is independent of $n$. The decomposition \eqref{eq.record_decomposition} writes the waiting index as one plus a sum of record gaps. To control the logarithm of this sum, we assign weights according to each record's height and position.

Partition the heights above $a$ into intervals
$$
\Delta_\ell:=(2^\ell a,2^{\ell+1}a],\qquad \ell\geq0.
$$
Give the initial term one in \eqref{eq.record_decomposition} weight $w_0:=1/2$. For $y\in\mathcal{R}\cap\Delta_\ell$, set
$$
w_y:=2^{-\ell-r_{\Delta_\ell}(y)-3},\qquad x_y:=D(y)\wedge n.
$$
Within each interval, the ranks $r_{\Delta_\ell}(y)$ are consecutive integers starting at zero. Summing the two geometric series shows that
$$
w_0+\sum_{\substack{y\in\mathcal{R}\\y>a}}w_y\leq1.
$$
Let
$$
S:=1+\sum_{\substack{y\in\mathcal{R}\\y>a}}x_y.
$$
The record decomposition gives $\kappa(a)\wedge n\leq S$. The log-sum inequality, applied to the terms $1$ and $x_y$ with their assigned weights, gives
$$
S\log S\leq\log(1/w_0)+
\sum_{\substack{y\in\mathcal{R}\\y>a}}x_y\log(x_y/w_y).
$$
Now expand the logarithms on the right. Since $u\log u$ is increasing for $u\geq1$, we obtain
$$
\begin{aligned}
(\kappa(a)\wedge n)\log(\kappa(a)\wedge n)
&\leq\sum_{\substack{y\in\mathcal{R}\\y>a}}x_y\log x_y+\log2\\
&\quad+(\log2)\sum_{\ell\geq0}
\sum_{y\in\mathcal{R}\cap\Delta_\ell}
(\ell+r_{\Delta_\ell}(y)+3)x_y.
\end{aligned}
$$

Take expectations. The change of root identity \eqref{eq.waiting_rerooting} turns the first sum into an integral of $h_n(y)$. In the last sum, use $x_y\leq D(y)$ and apply \eqref{eq.record_position_bound} to the terms involving ranks. Another application of \eqref{eq.waiting_rerooting} gives
$$
h_n(a)\leq\int_a^\infty\frac{h_n(y)}{\sinh y}\,\mathrm{d}y
+C+C\sum_{\ell\geq0}(\ell+1)
\int_{2^\ell a}^{2^{\ell+1}a}\frac{m(y)}{\sinh y}\,\mathrm{d}y.
$$
The remaining sum is bounded by a constant times $m(a)$. Indeed,
$$
m'(y)=-\frac{m(y)}{\sinh y},\qquad
m(2^\ell a)-1=\frac{2}{e^{2^\ell a}-1}
\leq2^{-\ell}(m(a)-1).
$$
Since $m(y)\to1$ as $y\to\infty$, telescoping after subtracting this limit gives
$$
\begin{aligned}
\sum_{\ell\geq0}(\ell+1)
\bigl(m(2^\ell a)-m(2^{\ell+1}a)\bigr)
&=\sum_{\ell\geq0}\bigl(m(2^\ell a)-1\bigr)\\
&\leq2(m(a)-1).
\end{aligned}
$$
We have proved
\begin{equation}
h_n(a)\leq C m(a)+\int_a^\infty\frac{h_n(y)}{\sinh y}\,\mathrm{d}y.
\label{eq.entropy_recursion}
\end{equation}

Iterate \eqref{eq.entropy_recursion} as in the mean waiting estimate. Since $h_n\leq n\log n$ and the kernel $1/\sinh y$ is integrable on $[a,\infty)$, the remaining term after $k$ substitutions is bounded by
$$
\frac{n\log n}{k!}
\left(\int_a^\infty\frac{\mathrm{d}y}{\sinh y}\right)^k,
$$
which tends to zero. Summing the terms generated by $C m$ gives
$$
\begin{aligned}
h_n(a)
&\leq C m(a)+C\int_a^\infty
\exp\left(\int_a^y\frac{\mathrm{d}z}{\sinh z}\right)
\frac{m(y)}{\sinh y}\,\mathrm{d}y\\
&=C m(a)(1+\log m(a)).
\end{aligned}
$$
For the last equality, use
$$
\exp\left(\int_a^y\frac{\mathrm{d}z}{\sinh z}\right)
=\frac{m(a)}{m(y)},\qquad
\int_a^\infty\frac{\mathrm{d}y}{\sinh y}=\log m(a).
$$
The bound is independent of $n$. Monotone convergence as $n\to\infty$ proves \eqref{eq.waiting_entropy}.
\end{proof}

\subsection{Backward waiting and record counts}

\begin{lemma}
There is a constant $C$ such that
\begin{equation}
\mathbb{E}_{\mathbb{Q}}
\left[\overline{\kappa}_{\mathrm b}(a)(1+N([a,2a]))\right]
\leq
\begin{cases}
C\log(e/a)/a,&0<a\leq1,\\
C,&a\geq1.
\end{cases}
\label{eq.mixed_waiting}
\end{equation}
\end{lemma}

\begin{proof}
Set
$$
s:=a\wedge1,\qquad I:=(s/2,s],\qquad J:=(a,2a],\qquad
V:=N(I)+N(I,J).
$$
Every height in $J$ exceeds every height in $I$. Thus, for each parent record in $I$, all of its child records in $J$ qualify for the paired count.

We use parent records in $I$ to relate backward waiting to forward waiting. Since $y\leq a$ for $y\in I$, monotonicity gives $\overline{\kappa}_{\mathrm b}(a)\leq\overline{\kappa}_{\mathrm b}(y)$. Integrating this inequality, then changing the root at each parent, yields
$$
\begin{aligned}
&\mathbb{E}_{\mathbb{Q}}
[\overline{\kappa}_{\mathrm b}(a)(1+N(J))]
\int_I\frac{\mathrm{d}y}{\sinh y}\\
&\quad\leq\int_I
\frac{\mathbb{E}_{\mathbb{Q}}
[\overline{\kappa}_{\mathrm b}(y)(1+N(J))]}{\sinh y}\,\mathrm{d}y\\
&\quad=\mathbb{E}_{\mathbb{Q}}\sum_{y\in\mathcal{R}\cap I}
p_y\bigl(1+N_{p_y,g_y}(J)\bigr)\\
&\quad\leq\mathbb{E}_{\mathbb{Q}}[\kappa(s/2)V].
\end{aligned}
$$
For the equality, use \eqref{eq.marked_rerooting} with
$$
\Phi(y,H,\eta,e):=
\mathbf{1}_I(y)\kappa_{\mathrm b}(y,\eta)(1+N(J)),
$$
where both statistics on the right are computed in $H$. Averaging the sign gives $\overline{\kappa}_{\mathrm b}(y)$ in the integral, while \eqref{eq.source_index} turns the backward waiting index into $p_y$ in the record sum. For the last inequality, every parent in $I$ occurs before the first positive level at most $s/2$, so $p_y<\kappa(s/2)$. The sum of the remaining factors $1+N_{p_y,g_y}(J)$ is exactly $V$.

The deterministic integral on the left is bounded away from zero uniformly in $a$. Since $s\leq1$ and $\sinh y\leq(\sinh1)y$ on $(0,1]$,
$$
\int_{s/2}^s\frac{\mathrm{d}y}{\sinh y}
\geq\frac{\log2}{\sinh1}.
$$
It therefore remains to bound $\mathbb{E}_{\mathbb{Q}}[\kappa(s/2)V]$.

The bounds for single records and pairs of records give exponential tails for both summands in $V$, with constants independent of $a$. A union bound yields
$$
\mathbb{Q}(V>u)\leq C e^{-cu},\qquad u\geq0,
$$
for some $c>0$. Reducing $c$ if necessary gives $\mathbb{E}_{\mathbb{Q}}e^{cV}\leq C$. We combine this exponential moment with the logarithmic moment of $\kappa(s/2)$ by the pointwise inequality
$$
cxV\leq x\log x-x+e^{cV},\qquad x>0.
$$
Taking $x:=\kappa(s/2)$ and applying \eqref{eq.waiting_entropy}, we obtain
$$
\mathbb{E}_{\mathbb{Q}}[\kappa(s/2)V]
\leq C m(s/2)(1+\log m(s/2))
\leq C\frac{\log(e/s)}s.
$$
For $a\leq1$, we have $s=a$, which gives the first range of \eqref{eq.mixed_waiting}. For $a\geq1$, we have $s=1$, which gives the constant bound. Finally, \eqref{eq.record_intensity} assigns zero expected count to any fixed height. Thus $N(J)=N([a,2a])$ almost surely for this fixed $a$, as required.
\end{proof}

\subsection{Returning to the independent exponential time}

We first compute the distribution of an excursion's value at an independent exponential time, with the additional requirement that the excursion has not yet reached a barrier $w$. For $w>0$ and a Borel set $A\subset(0,w)$, the formula is
\begin{equation}
\mathbf n_{\mathrm{pos}}\left(
\int_0^{\zeta\wedge T_w}\frac12e^{-s/2}
\mathbf{1}_{\{e_s\in A\}}\,\mathrm{d}s\right)
=\frac12\int_A\frac{\sinh(w-x)}{\sinh w}\,\mathrm{d}x.
\label{eq.excursion_endpoint}
\end{equation}
Here $s$ is the elapsed time since the excursion began, and $x$ is its absolute value at the observation time. Stopping the integral at $T_w$ enforces the barrier condition.

To verify the formula, consider Brownian motion started at $z\in(0,w)$ and killed when it leaves $(0,w)$. Its resolvent density at discount rate $1/2$ is
$$
r_w(z,x):=\frac{2\sinh(z\wedge x)\sinh(w-(z\vee x))}{\sinh w}.
$$
This means that its expected discounted occupation time satisfies
$$
\mathbb{E}_z\int_0^{T_0\wedge T_w}e^{-s/2}
\mathbf{1}_{\{B_s\in A\}}\,\mathrm{d}s
=\int_A r_w(z,x)\,\mathrm{d}x.
$$
As before, passage to positive excursion measure divides by $2z$ and lets $z\downarrow0$. The density of the exponential observation time contributes an additional factor $1/2$. The density in \eqref{eq.excursion_endpoint} is therefore
$$
\lim_{z\downarrow0}\frac{r_w(z,x)}{4z}
=\frac12\frac{\sinh(w-x)}{\sinh w}.
$$
The same formula holds for the absolute value of a negative excursion. If the barrier is removed by letting $w\to\infty$, each sign contributes the density $\frac12e^{-x}\,\mathrm{d}x$.

\begin{proof}[Proof of Lemma~\ref{lem.entrance_estimate}]
Fix $0<r\leq1$ and $0<\varepsilon\leq1/2$, and put $k:=\nu(r)$. Thus
$$
|B_\tau^{(i)}|>r\quad\text{for }0\leq i<k,\qquad
|B_\tau^{(k)}|=U(r)\leq r.
$$
Among the paths with indices $0,\ldots,k$, choose the one whose last zero before $\tau$ occurs latest. Denote its index by $j$, and define
$$
g:=\sup\{0\leq s<\tau:B_s^{(j)}=0\},\qquad
\eta:=\operatorname{sgn}(B_\tau^{(j)}),\qquad
x:=|B_\tau^{(j)}|.
$$
The choice of $j$ is almost surely unique. A tie would mean that an excursion containing $\tau$ starts at a zero of another coordinate. For each pair of coordinates, excursion compensation and Lemma~\ref{lem.disjoint_local_times} exclude such starts, by the same argument used in the change of root proof. There are only countably many pairs, so all ties are excluded simultaneously.

Every other inspected path therefore has its last zero strictly before $g$ and remains nonzero throughout $[g,\tau]$. Its sign and local time are constant on that interval, so its values can be described by an affine function of the reference excursion.

We first obtain a restriction on the possible value of $j$. For every original index $i<j$, the first entrance property gives $|B_\tau^{(i)}|>r$. If this lower path moves towards zero as the reference absolute value increases, affine propagation gives
$$
|B_g^{(i)}|-x=|B_\tau^{(i)}|>r.
$$
Its zero level in the direction of the reference excursion therefore exceeds $r+x$. If it moves away from zero, its signed level does not qualify for the backward search with sign $\eta$. None of the $j$ lower coordinates with original indices $0,\ldots,j-1$ can therefore stop that search at threshold $r+x$. In the history rooted at $(j,g)$, this says
\begin{equation}
j<\kappa_{\mathrm b,j,g}(r+x,\eta).
\label{eq.source_count_constraint}
\end{equation}

We estimate the entrance event by summing over all possible reference indices and their excursions. In the sums below, $i\geq0$ denotes a candidate value of the reference index $j$. The observation time inside an excursion starting at $g$ has the form $g+s$, and its density factors as
$$
\frac12e^{-(g+s)/2}=e^{-g/2}\frac12e^{-s/2}.
$$
Under excursion compensation, the first factor combines with local time at $g$ to give the root measure \eqref{eq.zero_measure}. Index stationarity then gives the same law $\mathbb{Q}$ for every candidate reference index. The second factor remains inside the excursion integral. Keeping necessary conditions for the entrance event and omitting the others gives an upper bound. We treat separately the cases $j=k$ and $j<k$.

Suppose first that $j=k$. The entrance occurs in the reference coordinate itself, so $x=U(r)$. Keep the restriction \eqref{eq.source_count_constraint} and require the absolute endpoint to lie in $[(1-\varepsilon)r,r]$. Compensation gives
$$
\begin{aligned}
&\mathbb{P}\bigl(j=k,\ (1-\varepsilon)r\leq U(r)\leq r\bigr)\\
&\quad\leq\mathbb{E}_{\mathbb{Q}}
\sum_{i\geq0}\sum_{\eta\in\{-1,1\}}
\mathbf n_{\mathrm{pos}}\left(
\int_0^\zeta\frac12e^{-s/2}
\mathbf{1}_{\{(1-\varepsilon)r\leq e_s\leq r\}}
\mathbf{1}_{\{i<\kappa_{\mathrm b}(r+e_s,\eta)\}}\,\mathrm{d}s\right).
\end{aligned}
$$
For a fixed history, endpoint $x$, and sign $\eta$, the number of integers $i\geq0$ satisfying the last indicator is exactly $\kappa_{\mathrm b}(r+x,\eta)$. Use the endpoint density with no barrier, then sum the two signs. This gives
\begin{equation}
\begin{aligned}
&\mathbb{P}\bigl(j=k,\ (1-\varepsilon)r\leq U(r)\leq r\bigr)\\
&\quad\leq\int_{(1-\varepsilon)r}^r
e^{-x}\mathbb{E}_{\mathbb{Q}}\overline{\kappa}_{\mathrm b}(r+x)\,\mathrm{d}x\\
&\quad\leq\varepsilon r\,
\mathbb{E}_{\mathbb{Q}}\overline{\kappa}_{\mathrm b}(r)\\
&\quad\leq C\varepsilon\log(e/r).
\end{aligned}
\label{eq.entrance_last_layer}
\end{equation}
The second inequality uses monotonicity of the backward waiting index, and the last follows from \eqref{eq.mixed_waiting}.

Now suppose that $j<k$. Since the reference coordinate precedes the first entrance, its endpoint satisfies $x>r$. Compute the higher zero levels in the history rooted at $(j,g)$. All inspected higher coordinates remain nonzero up to $\tau$, so their affine expressions remain valid. A negative zero level gives an absolute endpoint larger than $x$, and therefore larger than $r$. It cannot produce the entrance. The entrance coordinate must consequently have a positive zero level, which we denote by $w$.

The path at this positive level has not reached zero during the reference excursion. Its affine absolute value is $w-|B^{(j)}|$, so
$$
\max_{g\leq s\leq\tau}|B_s^{(j)}|<w,\qquad
U(r)=w-x\leq r.
$$
Every preceding positive zero level gives an endpoint larger than $r$, and hence exceeds $x+r$. Since $w\leq x+r$, the level $w$ is smaller than every preceding positive level. It is therefore a record in the history at $(j,g)$.

After shifting the reference index to zero, we can bound this case by summing over all record heights $w$ that satisfy these necessary conditions. For fixed $w$, define the allowed set of reference endpoints by
$$
A_w:=[w-r,w-(1-\varepsilon)r]\cap(r,\infty).
$$
The entrance band and $x>r$ require $x\in A_w$. This set lies in $(0,w)$ and has length at most $\varepsilon r$. The reference excursion must also be observed before it first reaches $w$. Retaining these conditions and \eqref{eq.source_count_constraint}, compensation gives
$$
\begin{aligned}
&\mathbb{P}\bigl(j<k,\ (1-\varepsilon)r\leq U(r)\leq r\bigr)\\
&\quad\leq\mathbb{E}_{\mathbb{Q}}
\sum_{i\geq0}\sum_{\eta\in\{-1,1\}}\sum_{w\in\mathcal{R}}
\mathbf n_{\mathrm{pos}}\left(
\int_0^{\zeta\wedge T_w}\frac12e^{-s/2}
\mathbf{1}_{\{e_s\in A_w\}}
\mathbf{1}_{\{i<\kappa_{\mathrm b}(r+e_s,\eta)\}}\,\mathrm{d}s\right).
\end{aligned}
$$
Each record height $w$ is fixed by the root history before the excursion is integrated. We may therefore apply \eqref{eq.excursion_endpoint} for that barrier, sum the indicators over $i$, and average the two signs. On $A_w$, we have $r+x\geq w$ and $w-x\leq r$. Monotonicity of $\overline{\kappa}_{\mathrm b}$ and the length bound for $A_w$ give
\begin{equation}
\begin{aligned}
&\mathbb{P}\bigl(j<k,\ (1-\varepsilon)r\leq U(r)\leq r\bigr)\\
&\quad\leq\mathbb{E}_{\mathbb{Q}}\sum_{w\in\mathcal{R}}
\int_{A_w}\overline{\kappa}_{\mathrm b}(r+x)
\frac{\sinh(w-x)}{\sinh w}\,\mathrm{d}x\\
&\quad\leq\varepsilon r\sinh r\,
\mathbb{E}_{\mathbb{Q}}\sum_{\substack{w\in\mathcal{R}\\w\geq r}}
\frac{\overline{\kappa}_{\mathrm b}(w)}{\sinh w}.
\end{aligned}
\label{eq.entrance_record_sum}
\end{equation}
Only heights $w>r$ can contribute, since otherwise $A_w$ is empty.

The backward waiting index and the record heights in this expectation come from the same history. This is the product controlled by \eqref{eq.mixed_waiting}. To apply that estimate, divide the heights into intervals $[b,2b)$ with $b:=2^\ell r$ and $\ell\geq0$. On each interval, both $\overline{\kappa}_{\mathrm b}(w)$ and $1/\sinh w$ are at most their values at $b$. Therefore
$$
\begin{aligned}
\mathbb{E}_{\mathbb{Q}}\sum_{w\in\mathcal{R}\cap[b,2b)}
\frac{\overline{\kappa}_{\mathrm b}(w)}{\sinh w}
&\leq\frac{\mathbb{E}_{\mathbb{Q}}
[\overline{\kappa}_{\mathrm b}(b)N([b,2b])]}{\sinh b}\\
&\leq
\begin{cases}
C\log(e/b)/b^2,&b\leq1,\\
C/\sinh b,&b>1.
\end{cases}
\end{aligned}
$$
Since $r\leq1$, we have $r\sinh r\leq C r^2$. The intervals with $b\leq1$ contribute to \eqref{eq.entrance_record_sum} at most
$$
C\varepsilon r^2\sum_{2^\ell r\leq1}
\frac{\log(e/(2^\ell r))}{(2^\ell r)^2}
\leq C\varepsilon\log(e/r)\sum_{\ell\geq0}4^{-\ell}
\leq C\varepsilon\log(e/r).
$$
For $b>1$, the sum of $1/\sinh b$ over the successive doubling values of $b$ is bounded by a constant independent of $r$. These intervals contribute at most $C\varepsilon r^2$, which is also bounded by $C\varepsilon\log(e/r)$. Combining the two cases with \eqref{eq.entrance_last_layer} proves Lemma~\ref{lem.entrance_estimate}.
\end{proof}

\paragraph{AI disclosure.} The author used GPT-6 Astra for mathematical discussions and manuscript drafting. The author directed this process, verified and revised the arguments, and takes full responsibility for the content.

\clearpage
\bibliographystyle{alpha}
\bibliography{ref}

\end{document}